\documentclass[11pt]{article}
\usepackage{latexsym}
\usepackage{amssymb}
\usepackage{amsmath}
\usepackage{amsfonts}
\usepackage{graphics}
\def\Bbb{\mathbb}

\def\eea{\end{eqnarray*}}

\newtheorem{thm}{Theorem}[section]
\newtheorem{prop}[thm]{Proposition}

\newtheorem{lem}[thm]{Lemma}

\newenvironment{proof}{\medskip \noindent
{\bf Proof.}}{\hfill \rule{.5em}{1em}
\\}

\newenvironment{rmk}{\mbox{ }\\{\bf  Remark}\mbox{ }}{
\hfill $\Box$\mbox{}\bigskip}
\begin{document}
\sloppy
\title{Surgery and equivariant Yamabe invariant}

\author{Chanyoung Sung
\thanks{email address: cysung@kias.re.kr \
}}

\date{Korea Institute for Advanced Study\\
207-43 Cheongryangri 2-dong Dongdaemun-gu\\
Seoul 130-722 Korea }

\maketitle

\begin{abstract}
We consider the equivariant Yamabe problem, i.e. the Yamabe
problem on the space of $G$-invariant metrics for a compact Lie
group $G$. The $G$-Yamabe invariant is analogously defined as the
supremum of the constant scalar curvatures of unit volume
$G$-invariant metrics minimizing the total scalar curvature
functional in their $G$-invariant conformal subclasses. We prove a
formula about how the $G$-Yamabe invariant changes under the
surgery of codimension $3$ or more, and compute some $G$-Yamabe
invariants.
\end{abstract}

\section{Introduction}
By the well-known uniformization theorem, the geometry and
topology of compact orientable surfaces have the trichotomy
according to the Euler characteristic. The Gauss-Bonnet theorem
says that the Euler characteristic is basically the constant
scalar curvature of the unit volume. Along this line one can
consider the following higher dimensional generalization,
so-called \emph{Yamabe invariant}.

Let $M$ be a smooth compact connected $n$-manifold. In analogy to
the $2$-dimension, let's consider the normalized Einstein-Hilbert
functional $$Q(g)= \frac{\int_M s_{g}\ dV_{g}}{(\int_M
dV_{g})^{\frac{n-2}{n}}}$$ defined on the space of smooth
Riemannian metrics on $M$, where $s_{g}$ and $dV_{g}$ respectively
denote the scalar curvature and the volume element of ${g}$. The
denominator is appropriately chosen for the purpose of the scale
invariance. But it turns out that this functional is neither
bounded above nor bounded below. In higher dimensions one need to
note that there are metrics which are not conformally equivalent
to each other. A \emph{conformal class} on M is by definition a
collection of smooth Riemannian metrics on M of the form
$$[g]\equiv \{\psi g \mid \psi:M \rightarrow \Bbb R^+ \},$$ where
$g$ is a fixed Riemannian metric. In each conformal class $[g]$
the above functional is bounded below and the minimum, called the
\emph{Yamabe constant} of $(M,[g])$ and denoted by $Y(M,[g])$, is
realized by a so-called \emph{Yamabe metric} which has constant
scalar curvature. By Aubin's theorem \cite{aubin}, the Yamabe
constant of any conformal class on any $n$-manifold is always
bounded by that of the unit $n$-sphere $S^n(1) \subset \Bbb
R^{n+1}$, which is $\Lambda_n\equiv
n(n-1)(\textrm{vol}(S^n(1)))^{2/n}$. The \emph{Yamabe invariant}
of $M$, $Y(M)$, is then defined as the supremum of the Yamabe
constant over the set of all conformal classes on $M$. Note that
it is a differential-topological invariant of $M$ depending only
on the smooth structure of the manifold.

The computation of the Yamabe invariant has been making notable
progress, particularly in low dimensions, due to LeBrun
\cite{lb1,lb2,lb3,IL}, Bray and Neves \cite{Bray}, Perelman
\cite{perel}, Anderson \cite{ander}, and etc. But in higher
dimensions little is known and noteworthy theorems to this end are
the surgery theorems. By the celebrated theorem of Gromov and
Lawson \cite{GL}, also independently by Schoen and Yau \cite{sy},
the Yamabe invariant of any manifold obtained from the manifolds
of positive Yamabe invariant by a surgery of codimension $3$ or
more is also positive. Moreover we have
\begin{thm}[Kobayashi \cite{koba}, Petean and Yun \cite{PY}]\label{know}
Let $M_1,M_2$ be smooth compact manifolds of dimension $n\geq 3$.
Suppose that an $(n-q)$-dimensional smooth compact (possibly
disconnected) manifold $W$ embeds into both $M_1$ and $M_2$ with
trivial normal bundle. Assume $q \geq 3$. Let $M$ be any manifold
obtained by gluing $M_1$ and $M_2$ along $W$. Then
$$
Y(M) \geq
\left\{
\begin{array}{ll}  -( |Y(M_1)|^{n/2}+  |Y(M_2)|^{n/2} )^{2/n}
   &\mbox{if } Y(M_i)\leq 0 \ \forall i\\
  \min(Y(M_1),Y(M_2)) &\mbox{if } Y(M_1) \cdot Y(M_2) \leq 0\\
  \min(Y(M_1),Y(M_2)) &\mbox{if } Y(M_i)\geq 0 \ \forall i\mbox{ and } q=n
\end{array}\right.
$$
\end{thm}
When $Y(M_i)\geq 0$ and $3\leq q\leq n-1$, no estimate has been
given even for $W=S^{n-q}$.

Now let's generalize this discussion to the equivariant Yamabe
problem. Let $G$ be a compact Lie group acting on $(M,g)$ smoothly
as an isometry. We will call such $(M,g)$ as a \emph{Riemannian}
$G$\emph{-manifold} and $[g]_G$ will denote the set of smooth
$G$-invariant metrics conformal to $g$. Then we have
\begin{thm}[Hebey and Vaugon \cite{HV}]\label{hebey}
Let $(M,g)$ a smooth compact Riemannian $G$-manifold. Then there
exists a metric $g'\in [g]_G$ of constant scalar curvature
realizing
$$Y(M,[g]_G):=\inf_{\hat{g}\in [g]_G}\frac{\int_M s_{\hat{g}}\
dV_{\hat{g}}}{(\int_M dV_{\hat{g}})^{\frac{n-2}{n}}},$$ and
$$Y(M,[g]_G)\leq \Lambda_n (\inf_{x\in M}|Gx|)^{\frac{2}{n}}$$ where $|Gx|$
denotes the cardinality of the orbit of $x$.
\end{thm}
We will call $Y(M,[g]_G)$ the $G$\emph{-Yamabe constant} of
$(M,[g]_G)$ and such a metric $g'$ will be called as a
$G$\emph{-Yamabe metric}. Obviously $Y(M,[g]_G)\geq Y(M,[g])$ for
any $G$-invariant metric $g$. We also remark that any $G$-Yamabe
metric with the nonpositive $G$-Yamabe constant is actually a
Yamabe metric, and hence the $G$-Yamabe constant coincides with
the Yamabe constant, because the constant scalar curvature metric
is unique up to constant in such a conformal class. The
$G$\emph{-Yamabe invariant} $Y_G(M)$ of $M$ is also defined as the
supremum of all the $G$-Yamabe constants. Of course it is an
invariant of the $G$-manifold $M$. We will show that some standard
theorems about the Yamabe constant can be generalized to the
$G$-Yamabe constant and prove the following surgery theorem for
the $G$-Yamabe invariant.
\begin{thm} Let $M_1,M_2$ be smooth
compact manifolds of dimension $n\geq 3$ on which a compact Lie
group $G$ acts smoothly.
Suppose that an $(n-q)$-dimensional smooth compact (possibly
disconnected) manifold $W$ with a locally transitive $G$-action
embeds $G$-equivariantly into both $M_1$ and $M_2$ with an
equivariant $G$-action on the trivial normal bundle. Assume $q
\geq 3$. Let $M$ be any $G$-manifold obtained by equivariantly
gluing $M_1$ and $M_2$ along $W$. Then
$$
Y_G(M) \geq \left\{
\begin{array}{ll}  -( |Y_G(M_1)|^{n/2}+  |Y_G(M_2)|^{n/2} )^{2/n}
   &\mbox{if } Y_G(M_i)\leq 0\ \forall i\\
  \min(Y_G(M_1),Y_G(M_2))
   &\mbox{otherwise. }
\end{array}\right.
$$
\end{thm}
In the final section we will use this to compute some $G$-Yamabe
invariants.
\section{Approximation of metric for Yamabe invariant}
Let's briefly go over the standard setup for the Yamabe problem.
Let $p=\frac{2n}{n-2}$, $a=4\frac{n-1}{n-2}$. Then
$$Q(\varphi^{p-2}g)=\frac{\int_M (a|d\varphi|_g^2 + s_g
\varphi^2)\ dV_g}{(\int_M |\varphi|^p dV_g)^{\frac{2}{p}}  },$$
and $$Y(M,[g]_G)=\inf \{Q(\varphi^{p-2}g) \mid \varphi \in
L_1^2(M) \ \textrm{is nonzero and } \ G\textrm{-invariant} \},$$
where the Sobolev space $L_1^2(M)$ is the set of $u\in L^2(M)$
such that $du \in L^2(M)$. A smooth function $\psi$ such that
$\psi^{p-2}g$ is a $G$-Yamabe metric will be called a
$G$\emph{-Yamabe minimizer} for $[g]_G$. Generalizing B. Bergery's
theorem \cite{ber}, the $G$-Yamabe constant also behaves
continuously with respect to the conformal class.
\begin{thm}\label{cont}
Let $g_i,g$ be $G$-invariant Riemannian metrics on $M$ such that
$g_i \rightarrow g$ in the $C^1$-topology, and $s_{g_i}
\rightarrow s_g$ in the $C^0$-topology on $M$. Then $Y(M,[g_i]_G)
\rightarrow Y(M,[g]_G)$.
\end{thm}
\begin{proof}
By the theorem \ref{hebey}, there exists a $G$-invariant conformal
change $\varphi^{p-2}g$ of $g$ making the scalar curvature
constant. Since $\varphi^{p-2}g_i \rightarrow \varphi^{p-2}g$ and
$s_{\varphi^{p-2}g_i}\rightarrow s_{\varphi^{p-2}g}$ for any
positive smooth function $\varphi$, we may assume that $s_g$ is
constant. We have two cases either $s_g\geq 0$, or $s_g< 0$.

Let's consider the first case. Given a sufficiently small
$\epsilon>0$, we can take an integer $N(\epsilon)$ such that for
$i\geq N(\epsilon),$ $$(1-\epsilon)g^{-1}\leq g^{-1}_i\leq
(1+\epsilon)g^{-1},$$
$$(1-\epsilon)dV_{g}\leq dV_{g_i}\leq (1+\epsilon)dV_g,$$ and
$$|s_{g}-s_{g_i}| \leq \epsilon.$$
Then for any $\varphi \in L_1^2(M)$
\begin{eqnarray*}
Q(\varphi^{p-2}g_i)&\leq& \frac{\int_M
(a(1+\epsilon)|d\varphi|_g^2 + (s_g+\epsilon)
\varphi^2)(1+\epsilon)\ dV_g}{(\int_M
|\varphi|^p(1-\epsilon)\ dV_g)^{\frac{2}{p}}}\\
&\leq&  \frac{(1+\epsilon)\int_M (a|d\varphi|_g^2 + s_g
\varphi^2)\ dV_g}{(1-\epsilon)^{\frac{2}{p}}(\int_M |\varphi|^p
dV_g)^{\frac{2}{p}}}+  \frac{\epsilon(1+\epsilon)\int_M
(a|d\varphi|_g^2 + \varphi^2)\
dV_g}{(1-\epsilon)^{\frac{2}{p}}(\int_M |\varphi|^p
dV_g)^{\frac{2}{p}}}\\ &\leq&
\frac{(1+\epsilon)}{(1-\epsilon)^{\frac{2}{p}}}Q(\varphi^{p-2}g)
+\frac{\epsilon(1+\epsilon)\bar{C}}{(1-\epsilon)^{\frac{2}{p}}},
\end{eqnarray*}
where $\bar{C}>0$ is a constant satisfying $$\int_M (a|d\psi|_g^2
+ \psi^2)\ dV_g \leq \bar{C}(\int_M |\psi|^p dV_g)^{\frac{2}{p}}$$
for any $\psi \in L_1^2(M)$, and similarly
\begin{eqnarray*}
Q(\varphi^{p-2}g_i)&\geq& \frac{\int_M
(a(1-\epsilon)|d\varphi|_g^2 + (s_g-\epsilon) \varphi^2)\
dV_{g_i}}{(\int_M |\varphi|^p\ dV_{g_i})^{\frac{2}{p}}}\\
&\geq& \frac{(1-\epsilon)\int_M (a|d\varphi|_g^2 + s_g \varphi^2)\
dV_g}{(1+\epsilon)^{\frac{2}{p}}(\int_M |\varphi|^p
dV_g)^{\frac{2}{p}}}-  \frac{\epsilon(1+\epsilon)\int_M
(a|d\varphi|_g^2 + \varphi^2)\
dV_g}{(1-\epsilon)^{\frac{2}{p}}(\int_M |\varphi|^p
dV_g)^{\frac{2}{p}}}\\ &\geq&
\frac{(1-\epsilon)}{(1+\epsilon)^{\frac{2}{p}}}Q(\varphi^{p-2}g)
-\frac{\epsilon(1+\epsilon)\bar{C}}{(1-\epsilon)^{\frac{2}{p}}}.
\end{eqnarray*}
Taking the infimum over $\varphi$ and letting $\epsilon
\rightarrow 0$, we get $Y(M,[g_i]_G)\rightarrow Y(M,[g]_G)$.

In the second case, we have $s_{g_i}<0$ for all sufficiently large
$i$. Recall O. Kobayashi's lemma \cite{koba}:
\begin{lem}\label{neg}
Let $(M,h)$ be any Riemannian $G$-manifold with $Y(M,[h]_G) \leq
0$. Then $$(\min s_h) \textrm{vol}_h(M)^{\frac{2}{n}} \leq
Y(M,[h]_G) \leq (\max s_h) \textrm{vol}_h(M)^{\frac{2}{n}}.$$
\end{lem}
\begin{proof}
The proof should be the same as the non-equivariant case because
$Y(M,[h]_G)=Y(M,[h])$ in this case. The case of $n=2$ is immediate
from the Gauss-Bonnet theorem. Let's consider the case when $n\geq
3$. The right inequality is obvious from
$$\frac{\int_M s_{h}\ dV_{h}}{(\int_M
dV_{h})^{\frac{n-2}{n}}}\leq (\max s_h)
\textrm{vol}_h(M)^{\frac{2}{n}}.$$ For the left inequality, we
claim that $\min s_h\leq 0$. Otherwise the Sobolev inequality says
that there exists a constant $\check{C}>0$ such that $(\int_M
\psi^p \ dV_h)^{\frac{2}{p}}\leq
\check{C}\int_M(a|d\psi|_h^2+s_h\psi^2) dV_h$ for any $\psi\in
L^2_1(M)$. This implies $Y(M,[h]_G) > 0$ which is contradictory to
the assumption. Once we have $\min s_h\leq 0$, by using the
H\"{o}lder inequality we get
$$(\min s_h) \textrm{vol}_h(M)^{\frac{2}{n}} \leq \frac{\int_M (\min
s_h)\varphi^2\ dV_h}{(\int_M |\varphi|^p dV_h)^{\frac{2}{p}} }\leq
Q(\varphi^{p-2}h)$$  for any $\varphi \in L_1^2(M)$, implying that
$(\min s_h) \textrm{vol}_h(M)^{\frac{2}{n}} \leq Y(M,[h]_G)$.
\end{proof}
By the above lemma, $$(\min s_{g_i})
\textrm{vol}_{g_i}(M)^{\frac{2}{n}} \leq Y(M,[g_i]_G) \leq (\max
s_{g_i}) \textrm{vol}_{g_i}(M)^{\frac{2}{n}}$$ for sufficiently
large $i$. Letting $i\rightarrow \infty$, we get
$Y(M,[g_i]_G)\rightarrow (\max s_{g})
\textrm{vol}_{g}(M)^{\frac{2}{n}}=Y(M,[g]_G)$.
\end{proof}

In the light of this, we want to find a sequence of $G$-invariant
metrics which has a nice form to perform a surgery and converges
to the given one. Generalizing the results of O. Kobayashi
\cite{koba}, and K. Akutagawa and B. Botvinnik \cite{aku}, we
present:
\begin{thm}
Let $W$ be a $G$-invariant submanifold of a Riemannian
$G$-manifold $(M,g)$  and let $\bar{g}$ be a $G$-invariant metric
defined in an open neighborhood of $W$, which coincides with $g$
on $W$ up to first derivatives, i.e. $g=\bar{g}$ and $\partial
g=\partial \bar{g}$ on $W$ and has the same scalar curvature as
$g$ on $W$. Then for sufficiently small $\delta
> 0$ there exists a $G$-invariant metric $g_\delta$ on $M$
satisfying the following properties.
\begin{description}
\item{(i)} $g_\delta\equiv g$ on $\{z\in M | \textrm{dist}_g(z,W)
> \delta \}$.
\item{(ii)} $g_\delta\equiv \bar{g}$ in an open neighborhood of
$W$.
\item{(iii)} $g_\delta \rightarrow g$ in the $C^1$-topology
on $M$ as $\delta \rightarrow 0$.
\item{(iv)} $s_{g_\delta}
\rightarrow s_g$ in the $C^0$-topology on $M$ as $\delta
\rightarrow 0$.
\end{description}
\end{thm}
\begin{proof}
Let $r$ be the $g$-distance from $W$. Obviously $r$ is
$G$-invariant. The proof goes in the same way as \cite{koba} and
\cite{aku}. We will be content with describing $g_\delta$. Given a
$\delta > 0$, take a smooth nonnegative function $w_\delta(r),
r\in [0, \infty)$ which satisfies $w_\delta(r)\equiv 1$ on
$[0,\frac{1}{4}e^{-\frac{1}{\delta}}]$, $w_\delta(r)\equiv 0$ on
$[\delta,\infty),$ $|r\frac{\partial w_\delta}{\partial r}| <
\delta,$ and $|r\frac{\partial^2 w_\delta}{\partial r^2}| <
\delta$. Then $g_\delta=g+w_\delta(r)(\bar{g}-g)$ does the job.
\end{proof}

To apply the above theorem we need to find a metric $\bar{g}$
which approximates $g$ near $W$ in a canonical way. Let's suppose
that $W$ has codimension $q$. Let
$(x,y)=(x^1,\cdots,x^{n-q},y^{n-q+1},\cdots,y^{n})$ be a local
trivialization of the normal bundle of $W$, where
$(x_1,\cdots,x_{n-q})$ is a local coordinate on the base $W$ and
$(y_{n-q+1},\cdots,y_q)$ is a coordinate on the fiber vector
space. Via the exponential map, this gives a local coordinate near
$W$. Let the indices $i,j,\cdots$ run from $1$ to $n-q$, and the
indices $\alpha,\beta,\gamma,\cdots$ run from $n-q+1$ to $n$.
Because we have taken the exponential normal coordinate in the
normal direction, we have on $W$
$$\frac{\partial}{\partial y^\alpha} g(\partial_i,\partial_j)=
g(\nabla_{\partial_\alpha}\partial_i,\partial_j)+
g(\partial_i,\nabla_{\partial_\alpha}\partial_j)=-2\Pi^\alpha_{ij},$$
\begin{eqnarray*}
\frac{\partial}{\partial y^\beta}
g(\partial_i,\partial_\alpha)&=&g(\nabla_{\partial_\beta}\partial_i,\partial_\alpha)+
g(\partial_i,\nabla_{\partial_\beta}\partial_\alpha)\\
&=&-g(\nabla_{\partial_i}\partial_\beta,\partial_\alpha)+g(\partial_i,0)=-\Gamma_{i\beta}^\alpha,
\end{eqnarray*}
and
\begin{eqnarray*}
\frac{\partial}{\partial y^\gamma}
g(\partial_\alpha,\partial_\beta)
=0,
\end{eqnarray*}
where
$\Pi^\alpha_{ij}=g(\partial_i,\nabla_{\partial_j}\partial_\alpha)$
is the second fundamental form of $W$, and
$\Gamma_{i\beta}^\alpha(x)$ is the Christoffel symbol for the
$g$-connection of the normal bundle on $W$. Therefore near $W$,
$g$ can be written as
\begin{align*}
g(x,y)=&\sum_{i,j}(g^W_{ij}(x)-2\sum_{\alpha}y^\alpha \Pi^\alpha_{ij}(x)+O(r^2))dx^idx^j\\
& + \sum_{i,\alpha,\beta}(-\Gamma_{i\beta}^\alpha(x)
y^\beta+O(r^2))dx^idy^\alpha+\sum_\alpha dy^\alpha
dy^\alpha+\sum_{\alpha \neq \beta}O(r^2)dy^\alpha dy^\beta,
\end{align*}
where $g^W=g|_W$ and $r=\sum_\alpha (y^\alpha)^2$. We will call
the above the canonical coordinate expression of $g$ near $W$.

Let $\hat{g}$ be the first order approximation of $g$, i.e.
$$\hat{g}:=\sum_{i,j}(g^W_{ij}(x)-2\sum_{\alpha}y^\alpha
\Pi^\alpha_{ij}(x))dx^idx^j+
\sum_{i,\alpha,\beta}(-\Gamma_{i\beta}^\alpha(x)
y^\beta)dx^idy^\alpha+\sum_\alpha dy^\alpha dy^\alpha.$$ Since $g$
and $r$ are $G$-invariant, $\hat{g}$ is also $G$-invariant. The
scalar curvature of $\hat{g}$ is in general different from that of
$g$. For the scalar curvature correction, we want to make a
conformal change which is $1$ at $W$ up to the first order. Let
$\bar{g}(x,y)=u(x,y)^{p-2}\hat{g}$ where $u$ is $G$-invariant,
\begin{eqnarray}\label{u}
u(x,0)=1,\ \ \ \textrm{and}\ \ \ \frac{\partial}{\partial
y^\alpha}u(x,0)=0
\end{eqnarray}
for any $\alpha$ on $W$.  Letting the uppercase Roman indices
denote $1$ through $n$ and using (\ref{u}), we have on $W$
\begin{eqnarray*}
\Delta_{\hat{g}} u&=&-\hat{\nabla}^A\partial_A
u=-\hat{g}^{AB}(\partial_A\partial_B
u-\hat{\Gamma}^C_{AB}\partial_C u)\\
&=&-\hat{g}^{\alpha\beta}\partial_\alpha\partial_\beta
u=-\sum_\alpha \frac{\partial }{\partial y^\alpha} \frac{\partial
u}{\partial y^\alpha},
\end{eqnarray*}
where $\hat{\nabla}$ and $\hat{\Gamma}$ denote the covariant
derivative and Christoffel symbol of $\hat{g}$ respectively. We
set
$$u(x,y):=1-\frac{r^2}{8aq}(s_g|_W-s_{\hat{g}}|_W).$$ Then
on $W$,
$$s_{\bar{g}}=u^{1-p}(4a\Delta_{\hat{g}} u+s_{\hat{g}}u)
=-4a\sum_\alpha \frac{\partial }{\partial y^\alpha} \frac{\partial
u}{\partial y^\alpha}+s_{\hat{g}}=s_g.$$ Combined with the above
theorem, we obtain :
\begin{thm}\label{g-delta}
Let $W$ be a $G$-invariant submanifold of a Riemannain
$G$-manifold $(M,g)$. For sufficiently small $\delta>0$, there
exists a $G$-invariant metric $g_\delta$ such that
\begin{description}
\item{(i)} $g_\delta \rightarrow g$ in the $C^1$-topology on $M$
as $\delta \rightarrow 0$.
\item{(ii)} $s_{g_\delta} \rightarrow
s_g$ in the $C^0$-topology on $M$ as $\delta \rightarrow 0$.
\item{(iii)} $g_\delta\equiv g$ on $\{z\in M |
\textrm{dist}_g(z,W)
> \delta \}$.
\item{(iv)} In an open neighborhood of $W$, $g_\delta$ is
conformally equivalent to
$\sum_{i,j}(g^W_{ij}(x)-2\sum_{\alpha}y^\alpha
\Pi^\alpha_{ij}(x))dx^idx^j+
\sum_{i,\alpha,\beta}(-\Gamma_{i\beta}^\alpha(x)
y^\beta)dx^idy^\alpha+\sum_\alpha dy^\alpha dy^\alpha$.
\end{description}

\end{thm}

For the conformal classes which are close in a $G$-invariant
subset, we can obtain a common upper bound.
\begin{prop}\label{bound}
Let $\{g_\alpha|\alpha\in I \}$ be a collection of smooth
$G$-invariant metrics on a compact $G$-manifold $X$. Suppose that
there exists a constant $D_1$ and $D_2$ such that
$|g_\alpha-g_\beta| \leq D_1$ and $|s_{g_\alpha}-s_{g_\beta}| \leq
D_2$ in some $G$-invariant open subset $U\subset X$ for any
$\alpha, \beta \in I$. Then there exists a constant $D$ such that
$Y(X,[g_\alpha]_G)\leq D$ for any $\alpha\in I$.
\end{prop}
\begin{proof}
Take a smooth bump function $\phi(x)\geq 0$ supported in $U$. In
general $\phi$ is not $G$-invariant. Let $d\mu$ be the unit-volume
bi-invariant measure on $G$. Define $\bar{\phi}(x):=\int_G
\phi(gx)\ d\mu(g).$ Then $\bar{\phi}$ is $G$-invariant and also
supported in $U$. Now $Q(\bar{\phi}^{p-2}g_\alpha)$ is bounded
above and by definition $Y(X,[g_\alpha]_G)\leq
Q(\bar{\phi}^{p-2}g_\alpha)$ for any $\alpha\in I$.

\end{proof}

\section{Proof of Main theorem}
We start with the equivariant version of O. Kobayashi's lemma
\cite{koba}.
\begin{lem}
Let $(M_1 \cup M_2,g_1 \cup g_2)$ be the disjoint union of
$(M_1,g_1)$ and $(M_2,g_2)$. Then $Y(M_1\cup  M_2,[g_1 \cup
g_2]_G)$ is given by
$$\left\{
\begin{array}{ll}  -( |Y(M_1,[g_1]_G)|^{n/2}+  |Y(M_2,[g_2]_G)|^{n/2} )^{2/n}
   &\mbox{if } Y(M_i,[g_i]_G)\leq 0\ \forall i\\
  \min(Y(M_1,[g_2]_G),Y(M_2,[g_2]_G)) &\mbox{otherwise,}
\end{array}\right.
$$
and
$$
Y_G(M_1 \cup M_2) = \left\{
\begin{array}{ll}  -( |Y_G(M_1)|^{n/2}+  |Y_G(M_2)|^{n/2} )^{2/n}
   &\mbox{if } Y_G(M_i)\leq 0\ \forall i\\
  \min(Y_G(M_1),Y_G(M_2)) &\mbox{otherwise. }
\end{array}\right.
$$
\end{lem}
\begin{proof}
Suppose $Y(M_1,[g_1]_G)\geq Y(M_2,[g_2]_G)\geq 0$. Then for any
$c^2g_1' \cup g_2'\in [g_1 \cup g_2]_G$ where $c>0$ is a constant,
\begin{align*}
Q(c^2g_1' \cup g_2')&= \frac{\int_{M_1} c^{n-2}s_{g_1'}\ dV_{g_1'}
+\int_{M_2}s_{g_2'}\ dV_{g_2'} } {(\int_{M_1}c^n\
dV_{g_1'}+\int_{M_2} dV_{g_2'})^{\frac{n-2}{n}} }\\  &\geq
\frac{(\int_{M_1}c^n\
dV_{g_1'})^{\frac{n-2}{n}}Y(M_1,[g_1]_G)+(\int_{M_2}
dV_{g_2'})^{\frac{n-2}{n}}Y(M_2,[g_2]_G)}{(\int_{M_1}c^n\
dV_{g_1'}+\int_{M_2} dV_{g_2'})^{\frac{n-2}{n}} }\\ &\geq
\frac{(\int_{M_1}c^n\
dV_{g_1'})^{\frac{n-2}{n}}Y(M_2,[g_2]_G)+(\int_{M_2}
dV_{g_2'})^{\frac{n-2}{n}}Y(M_2,[g_2]_G)}{(\int_{M_1}c^n\
dV_{g_1'})^{\frac{n-2}{n}}+(\int_{M_2} dV_{g_2'})^{\frac{n-2}{n}}}\\
&= Y(M_2,[g_2]_G),
\end{align*}
and $Q(c^2g_1' \cup g_2')\rightarrow Y(M_2,[g_2]_G)$ if
$c\rightarrow 0$ and $g_2'$ is a $G$-Yamabe metric on $M_2$.

Suppose $Y(M_1,[g_1]_G)\geq 0 \geq Y(M_2,[g_2]_G)$. Also for any
$c^2g_1' \cup g_2'\in [g_1 \cup g_2]_G$,
\begin{eqnarray*}
Q(c^2g_1' \cup g_2')&=& \frac{\int_{M_1} c^{n-2}s_{g_1'}\
dV_{g_1'} +\int_{M_2} s_{g_2'}\ dV_{g_2'}} {(\int_{M_1}c^n\
dV_{g_1'}+\int_{M_2} dV_{g_2'})^{\frac{n-2}{n}}} \\ &\geq&
\frac{0+(\int_{M_2}
dV_{g_2'})^{\frac{n-2}{n}}Y(M_2,[g_2]_G)}{(\int_{M_1}c^n\
dV_{g_1'}+\int_{M_2}
dV_{g_2'})^{\frac{n-2}{n}}}\\
&\geq& \frac{(\int_{M_2}
dV_{g_2'})^{\frac{n-2}{n}}Y(M_2,[g_2]_G)}{(\int_{M_2}dV_{g_2'})^{\frac{n-2}{n}}}\\
&=& Y(M_2,[g_2]_G),
\end{eqnarray*}
and $Q(c^2g_1' \cup g_2')\rightarrow Y(M_2,[g_2]_G)$ if
$c\rightarrow 0$ and $g_2'$ is a $G$-Yamabe metric on $M_2$.

For the last remaining case, suppose $Y(M_i,[g_i]_G)\leq 0$ and we
assume $g_i$ is a $G$-Yamabe metric for $(M_i,[g_i]_G)$ for each
$i$ such that $s_{g_1}=s_{g_2}<0$. Now note that the lemma
\ref{neg} still holds true for the non-connected manifolds and its
corollary is that any $G$-invariant metric of nonpositive constant
scalar curvature is a $G$-Yamabe metric. Thus $g_1\cup g_2$ is a
Yamabe metric and
\begin{eqnarray*}
Y(M_1 \cup M_2, [g_1\cup g_2]_G)&=&s_{g_1\cup g_2}\
\textrm{vol}_{g_1\cup g_2}(M_1 \cup M_2)^{\frac{2}{n}}\\ &=&
-(|s_{g_1\cup g_2}|^{\frac{n}{2}} \textrm{vol}_{g_1\cup g_2}(M_1
\cup M_2))^{\frac{2}{n}}\\ &=& -(|s_{g_1}|^{\frac{n}{2}}
\textrm{vol}_{g_1}(M_1)+|s_{g_2}|^{\frac{n}{2}}
\textrm{vol}_{g_2}(M_2))^{\frac{2}{n}}\\
&=&-(|Y(M_1,[g_1]_G)|^{\frac{n}{2}}+
|Y(M_2,[g_2]_G)|^{\frac{n}{2}})^{\frac{2}{n}}.
\end{eqnarray*}

The second assertion is immediately obtained by taking the
supremum of the first equality.
\end{proof}

By the above lemma, we only need to prove the following theorem.
\begin{thm}
Let $M_0$ be a smooth compact (possibly disconnected) manifold of
dimension $n\geq 3$ on which a compact Lie group $G$ acts
smoothly,
and $W$ be an $(n-q)$-dimensional smooth compact (possibly
disconnected) manifold with a locally transitive $G$-action.
Suppose that two copies of $W$ embed $G$-equivariantly into $M_0$
with an equivariant $G$-action on the trivial normal bundle.
Assume $q \geq 3$. Let $M$ be any $G$-manifold obtained by an
equivariant surgery on $M_0$ along $W$. Then
$$
Y_G(M) \geq Y_G(M_0).
$$
\end{thm}
\begin{proof}
The idea of proof when $q=n$ is the same as the well-known result
of Osamu Kobayashi \cite{koba}, which considers a gluing with a
long neck. When $q<n$, the idea is inspired by Dominic Joyce's
method in \cite{joyce}. We construct $M$ with the volume of the
gluing region very small. This forces the $G$-Yamabe minimizer of
$M$ to concentrate away from the gluing region, otherwise the
value of Yamabe functional gets too big. Then the $G$-Yamabe
constant of $M$ is basically expressed by that of $M_0$. Although
we can simplify our proof a little bit by restricting to the case
$Y_G(M_0)>0$, we will prove the general case for completeness.  By
abuse of notation $W$ will also denote the submanifolds embedded
in $M$.

Let $0< \epsilon_1, \epsilon_2\ll 1$. Take a conformal class
$[g_0]_G$ on $M_0$ such that $Y(M_0,[g_0]_G)\geq Y_G(M_0)-
\frac{\epsilon_1}{2}$. Applying the theorem \ref{cont} and
\ref{g-delta}, we can find a $G$-invariant metric $g$ satisfying
$Y(M_0,[g]_G)\geq Y(M_0,[g_0]_G)- \frac{\epsilon_1}{2}$ and $g$
near $W$ is the canonical first order approximation of $g_0$, i.e.
$g=\sum_{i,j}((g^W_0)_{ij}(x)-2\sum_{\alpha}y^\alpha
\Pi^\alpha_{ij}(x))dx^idx^j+
\sum_{i,\alpha,\beta}(-\Gamma_{i\beta}^\alpha(x)
y^\beta)dx^idy^\alpha+\sum_\alpha dy^\alpha dy^\alpha$,where
$(y_{n-q+1},\cdots,y_n)$ is the $g_0$-exponential normal
coordinate in the normal direction. Since
$\textrm{dist}_{g}((x,y),W)=\sum_\alpha(y^\alpha)^2$, it turns out
that $(y_{n-q+1},\cdots,y_n)$ is also the $g$-exponential normal
coordinate, and so the above expression of $g$ is the canonical
coordinate expression for $g$ itself by the uniqueness. So we may
assume that
$$Y(M_0,[g]_G)\geq Y_G(M_0)- \epsilon_1$$ and
$$g=\sum_{i,j}(g^W_{ij}(x)-2\sum_{\alpha}y^\alpha
\Pi^\alpha_{ij}(x))dx^idx^j+
\sum_{i,\alpha,\beta}(-\Gamma_{i\beta}^\alpha(x)
y^\beta)dx^idy^\alpha+\sum_\alpha dy^\alpha dy^\alpha$$ on
$N(r_0):= \{ r=(\sum_\alpha y_\alpha^2)^{\frac{1}{2}} \leq r_0
\}$. Also keep in mind that the $G$-action fixes $r$, and acts on
$x$ as in $W$.

We first consider the case when $q=n$, i.e. $W$ is a finite set of
points. In this case $g$ is the Euclidean metric near $W$. Since
$r$ is $G$-invariant, by multiplying a conformal factor $f(r)$
which is $\frac{1}{r^2}$ near $W$, $(M_0-(W\cup W),g)$ is
conformal to a Riemannian $G$-manifold $(M_0',g')$ whose end is
two copies of an infinite cylinder $W\times S^{n-1}(1) \times
[0,\infty)$. Cut off both infinite cylinders at a large integer
$l\in [0,\infty)$ and glue them along the boundary to get a
Riemannian $G$-manifold $(M_l,\bar{g}_l)$ which contains a
cylinder $W\times S^{n-1}(1) \times [0,2l]$. Note that the
complement of the cylindrical region in $M_l$ is $G$-invariant and
the same for any $l$. Thus by the proposition \ref{bound},
$\{Y(M_l,[\bar{g}_l]_G)| \ l\in [0,\infty)\}$ is bounded above.
This is an important fact to be used below.

To estimate a lower bound of $Y(M_l,[\bar{g}_l]_G)$, let $\psi_l$
be a $G$-Yamabe minimizer satisfying $\int_{M_l} \psi^p_l\
dV_{\bar{g}_l}=1.$ Since $\{Y(M_l,[\bar{g}_l]_G)| \ l\in
[0,\infty)\}$ is bounded above, there exists a constant $A>0$
independent of $l$ such that $$\int_{W\times S^{n-1}(1) \times
[0,2l]} (a|d\psi_l|_{\bar{g}_l}+2(n-1)(n-2)\psi_l^2)\
dV_{\bar{g}_l}\leq A.$$ Combined with $\int_{W\times S^{n-1}(1)
\times [0,2l]} \psi^p_l\ dV_{\bar{g}_l}<1$, it implies that there
exists an integer $N_l \in [0,l-1]$ such that
\begin{equation}\label{cyl1}
\int_{W\times S^{n-1}(1) \times [2N_l,2N_l+2] }
(a|d\psi_l|_{\bar{g}_l}+2(n-1)(n-2)\psi_l^2)\ dV_{\bar{g}_l}\leq
\frac{A+1}{l},
\end{equation}
and
\begin{equation}\label{cyl2}
\int_{W\times S^{n-1}(1) \times [2N_l,2N_l+2] } \psi_l^p\
dV_{\bar{g}_l}\leq \frac{A+1}{l}.
\end{equation}
Let $\xi(t):\Bbb R \rightarrow [0,1]$ be a smooth function such
that
$$\xi(t)=\left\{ \begin{array}{cl}
1& \mbox{for } t\in (-\infty,0]\cup [2,\infty)\\
0&\mbox{for } t\in [\frac{2}{3},\frac{4}{3}].
\end{array}\right.$$
Define a smooth function $\Psi_l$ on $M_l$ as
$$\Psi_l=\left\{ \begin{array}{cc}
\psi_l(z,t)\xi(t-2N_l)& \mbox{for } (z,t)\in (W\times
S^{n-1}(1))\times [0,2l]\\
\psi_l&\mbox{elsewhere. }
\end{array}\right.$$
Cut $M_l$ at $W\times S^{n-1}(1)\times \{2N_l+1\}$ and glue two
half infinite cylinders to get back $(M_0',g')$. Extend $\Psi_l$
to $M_0'$ by defining it to be zero on the additional half
infinite cylinders. Noting (\ref{cyl1}), (\ref{cyl2}), and the
fact that $\{Y(M_l,[\bar{g}_l]_G)| \ l\in [0,\infty)\}$ is bounded
above,
we can get $$Q(\Psi_l^{p-2}g')\leq
Y(M_l,[\bar{g}_l]_G)+\frac{B}{l},$$ where $B$ is a constant
independent of $l$. This implies $Y_G(M_0)- \epsilon_1\leq
Y(M_0,[g]_G)\leq Y(M_l,[\bar{g}_l]_G)+\frac{B}{l}.$ Letting
$l\rightarrow \infty$ and $\epsilon_1\rightarrow 0$, we finally
obtain $Y_G(M_0)\leq Y_G(M).$

Now we turn to the case of $q<n$ which will be needed at the last
stage. We will perform a refined version of the well-known
Gromov-Lawson bending \cite{GL,Sung1} on $N(r_0)$. The manifold is
constructed as a hypersurface in the Riemannian product $\Bbb R
\times M_0$ in accordance with an appropriate smooth curve
$\gamma$ in $\{(t,r)\in \Bbb R^2 \}$, which starts tangentially to
the $r$-axis at $t=0$ and ends up parallel to the $t$-axis as in
the following figure. We extend the isometric $G$-action to $\Bbb
R \times M_0$ in an obvious way that $t$ is invariant. Since $r$
is $G$-invariant, the constructed manifold is a $G$-invariant
submanifold of the Riemannian $G$-manifold, and hence also a
Riemannian $G$-manifold. The angle of bending at each radius is
denoted by $\theta$, and $k\geq 0$ denotes the geodesic curvature.
\begin{figure}[h]
  \centering
  \scalebox{0.9}{\includegraphics{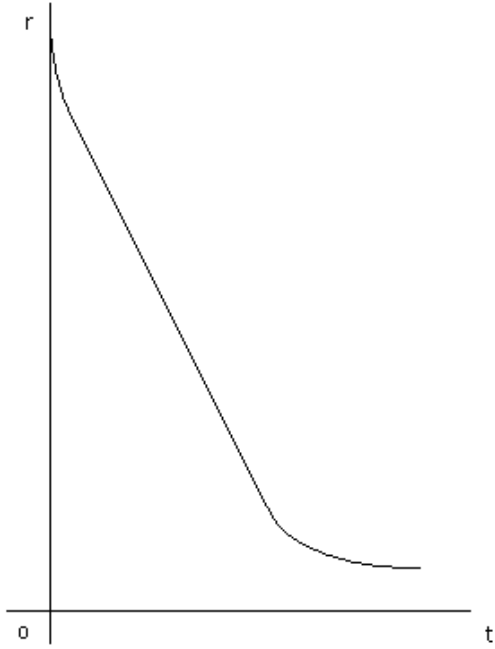}}
  \caption{curve $\gamma$}
\end{figure}
The scalar curvature $s$ is given by
\begin{eqnarray*}
s&=&s_g-2Ric_g(\frac{\partial}{\partial
r},\frac{\partial}{\partial
r})\sin^2\theta+(-\frac{2(q-1)}{r}+O(1))k\sin\theta\\& &
+(q-1)(q-2)\frac{\sin^2\theta}{r^2}+O(1)\frac{\sin^2\theta}{r}\\
&\geq&s_g+\frac{(q-1)(q-2)}{2}\frac{\sin^2\theta}{r^2}-3(q-1)\frac{k\sin\theta}{r},
\end{eqnarray*}
for sufficiently small $r>0$, where $s_g$ and $Ric_g$ denote the
scalar curvature and the Ricci curvature of $g$ respectively.

The construction of $\gamma$ is done in 3 steps. First, by
continuity we make a bending of small $\theta_0$ keeping
$$\frac{(q-1)(q-2)}{2}\frac{\sin^2\theta}{r^2}-3(q-1)\frac{k\sin\theta}{r}
> -\epsilon_2$$ so that $s> s_g-\epsilon_2$. Let $r_1$ be the
radius at the end and take $r_1'$ such that $0<r_1' \ll r_1$. As a
second step $\gamma$ goes down to $r=r_2$ straight i.e. $k=0$.
Since $k=0$, we have in this step $$s\geq
s_g+\frac{(q-1)(q-2)}{2}\frac{\sin^2\theta_0}{r^2}> s_g.$$ Here
$r_2>0$ is chosen small enough so that there exists a $C^\infty$
function $\eta(r):\Bbb R^+\rightarrow [0,1]$ such that
$$\eta(r)=\left\{ \begin{array}{cl}
0&\mbox{for } r\leq r_2\\
1& \mbox{for } r\geq r_1',
\end{array}\right.$$
and $$|d\eta| \leq \sqrt{\frac{(q-1)(q-2)}{2}}\frac{\sin
\theta_0}{r}.$$ (Consider the graph of
$y=(\sqrt{\frac{(q-1)(q-2)}{2}}\sin \theta_0)\ln x$.) This
$\eta(r)$ will be used later as a radial cut-off function on
$(M_0,g)$. Now the third step proceeds. We bend $\gamma$ after the
following prescription of the curvature function $k(L)$
parameterized by the arc length $L$.
\begin{figure}[h]
  \centering
  \scalebox{0.8}{\includegraphics{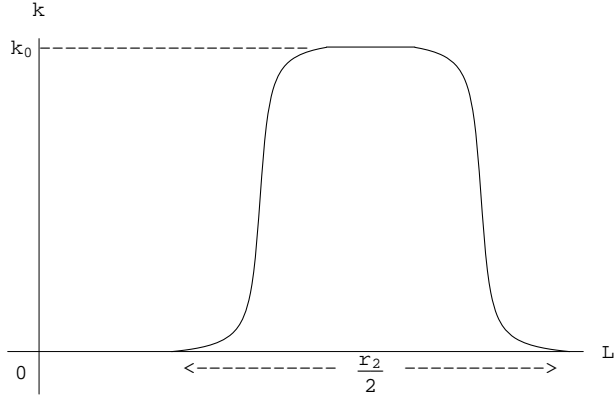}}
  \caption{curvature function $k(L)$}
\end{figure}
Here, $k_0$, the maximum of $k$, is defined as
$\frac{(q-2)\sin\theta_0}{6r_2}$ so that
\begin{equation}\label{original}
\frac{(q-1)(q-2)}{2} \frac{\sin^2\theta}{r^2}
-3(q-1)\frac{k\sin\theta}{r} \geq 0
\end{equation}
is ensured during this process and hence $s \geq s_g $. The amount
of the bend $\Delta\theta$ is
$$\Delta\theta=\int k \ dL\approx
k_0\cdot\frac{r_2}{2}=\frac{(q-2)\sin\theta_0}{12}.$$ Repeat this
process with the curvature prescription completely determined only
by the ending radius of the previous process until we achieve a
total bend of $\frac{\pi}{2}$. So the length of $\gamma$ during
this step is less than
\begin{eqnarray}
\frac{r_2}{2}([\frac{\pi}{2}/ \Delta\theta]+1) \leq \frac{3\pi
r_2}{(q-2)\sin \theta_0}+\frac{r_2}{2}.\label{leng}
\end{eqnarray}
Let $r_3$ be the final radius.

To smoothly glue two bent regions along the boundary $W \times
S^{q-1}$, we have to homotope the metrics on the boundaries. Let
$h_{r}$ be the metric on $W \times S^{q-1}$ induced from the
boundary of $(N(r),g)$. On $W \times S^{q-1}$ we define a
$G$-invariant product metric
$\bar{h}_{r}:=\sum_{i,j}\bar{g}^W+g_{std}(r)$ where $\bar{g}^W$ is
a fixed $G$-invariant metric on $W$ and $g_{std}(r)$ denotes the
round metric of $S^{q-1}(r)$. Obviously the scalar curvature
$s_{\bar{h}_r}$ of $\bar{h}_r$ is $\frac{(q-1)(q-2)}{r^2}+O(1)$.
Moreover
\begin{lem}\label{zz}
Let $h_r^\nu$ for $\nu\in [0,1]$ be the convex combination $\nu
h_r+(1-\nu)\bar{h}_{r}$ of $h_r$ and $\bar{h}_{r}$. Then there
exists a constant $C>0$ such that the scalar curvature
$s_{h_r^\nu}$ of $h_r^\nu$ is bounded below $\frac{C}{r^2}$ for
any $\nu$ and any sufficiently small $r>0$.
\end{lem}
\begin{proof}
This is basically because $h_r^\nu$ is very close to a riemannian
submersion with totally geodesic fibers $S^{q-1}(r)$, and hence
the O'Neill's formula \cite{bes} gives such an estimate of
$s_{h_r^\nu}$. It's enough to show that the difference between
$s_{h_r^\nu}$ and $s_{\bar{h}_r}$ is at most $O(\frac{1}{r})$.

As before we let $i,j,k,\cdots$ denote the indices of coordinates
of $W$ in $W\times S^{q-1}$ and $\alpha,\beta,\gamma,\cdots$
denote the indices of coordinates of $S^{q-1}$ in $W\times
S^{q-1}$, and $A,B,C,\cdots$ will denote the indices of
coordinates of both $W$ and $S^{q-1}$. Writing an $(n-1)\times
(n-1)$ matrix $(M_{AB})$ as
$$ \left(
\begin{tabular}{c|c}
&\\
$M_{ij}$&$M_{i\alpha}$\\ &\\
\cline{1-2}&\\
$M_{\alpha i}$ & $M_{\alpha\beta}$\\&\\
\end{tabular}
 \right),
 $$
 we have
\begin{eqnarray}\label{mat1}
(h_r^\nu)&=&(\bar{h}_r)+(h_r^\nu-\bar{h}_r) \nonumber\\ &=& \left(
\begin{tabular}{c|c}
&\\
$O(1)$&$0$\\ &\\
\cline{1-2}&\\
$0$ & $O(r^2)$\\&\\
\end{tabular}
 \right)
 +\left(
\begin{tabular}{c|c}
&\\
$O(r)$&$O(r^2)$\\ &\\
\cline{1-2}&\\
$O(r^2)$ & $0$\\&\\
\end{tabular}
 \right),
\end{eqnarray}
and
\begin{eqnarray}\label{mat2}
(h_r^\nu)^{-1}&=&(\bar{h}_r)^{-1}+((h_r^\nu)^{-1}-(\bar{h}_r)^{-1})\nonumber\\
&=& \left(
\begin{tabular}{c|c}
&\\
$O(1)$&$0$\\ &\\
\cline{1-2}&\\
$0$ & $O(\frac{1}{r^2})$\\&\\
\end{tabular}
 \right)
+ \left(
\begin{tabular}{c|c}
&\\
$O(r)$&$O(1)$\\ &\\
\cline{1-2}&\\
$O(1)$ & $O(\frac{1}{r})$\\&\\
\end{tabular}
 \right).
\end{eqnarray}
The same estimates also hold for their derivatives.

Recall that Christoffel symbols of a metric $h$ are given by
\begin{equation}\label{chris}
\Gamma_{AB}^C=\frac{1}{2}\sum_D h^{CD} \{ \frac{\partial
h_{AD}}{\partial x_B}+ \frac{\partial h_{BD}}{\partial x_A}-
\frac{\partial h_{AB}}{\partial x_D}\},
\end{equation}
and the Riemann curvature tensor $R$ is given by
\begin{equation}\label{riemann}
R_{ABC}^D=\partial_A\Gamma_{BC}^D-
\partial_B\Gamma_{AC}^D+
\Gamma_{BC}^E\Gamma_{AE}^D- \Gamma_{AC}^E\Gamma_{BE}^D.
\end{equation}
Denote the Christoffel symbol of $\bar{h}_r$ and $h_r^\nu$ by
$\bar{\Gamma}_r$ and $\Gamma_r^\nu$ respectively. Then the direct
computations show that
$$(\bar{\Gamma}_r)_{AB}^C=O(1)=\partial(\bar{\Gamma}_r)_{AB}^C,$$
and
$$(\Gamma_r^\nu)_{AB}^C-(\bar{\Gamma}_r)_{AB}^C=O(r)=
\partial(\Gamma_r^\nu)_{AB}^C-\partial(\bar{\Gamma}_r)_{AB}^C$$ except
$$(\Gamma_r^\nu)^\alpha_{ij}-(\bar{\Gamma}_r)^\alpha_{ij}=O(\frac{1}{r})
=\partial(\Gamma_r^\nu)^\alpha_{ij}-\partial(\bar{\Gamma}_r)^\alpha_{ij}.$$
Also denote the Riemann curvature tensor of $\bar{h}_r$ and
$h_r^\nu$ by $\bar{R}_r$ and $R_r^\nu$ respectively. Then
$$(R_r^\nu)^\alpha_{\alpha ij}-(\bar{R}_r)^\alpha_{\alpha
ij}=O(\frac{1}{r})=(R_r^\nu)^l_{ijk}-(\bar{R}_r)^l_{ijk},$$ and
$$(R_r^\nu)^\alpha_{\alpha\beta\gamma}-(\bar{R}_r)^\alpha_{\alpha\beta\gamma}=O(r).$$
Thus the difference between sectional curvatures of $\bar{h}_r$
and $h_r^\nu$ is bounded above by $O(\frac{1}{r})$, and hence so
is the differences of two scalar curvatures, completing the proof.
\end{proof}

Now we have the metric $h_{r_3}$ on the boundary. We have to
homotope $h_{r_3}$ to a $G$-invariant product metric
$\bar{h}_{r_3}$. Consider a smooth homotopy
$H_{r_3}(z,t):=\varphi(t)h_{r_3}+(1-\varphi(t)) \bar{h}_{r_3}$ for
$(z,t)\in (W\times S^{q-1}) \times [0,1]$, where
$\varphi:[0,1]\rightarrow [0,1]$ is a smooth decreasing function
which is $1$ near $0$ and $0$ near $1$. In the above lemma we have
seen that $(W\times S^{q-1},H_{r_3}(z,t))$ for each $t\in [0,1]$
has positive scalar curvature. Then by the Gromov-Lawson lemma in
\cite{GL}, there exists a constant $d>0$ such that the metric
$H_{r_3}(z,t/d)+dt^2$ on $W \times S^{q-1} \times [0,d]$ has
positive scalar curvature for sufficiently small $r_3>0$.
Obviously $H_{r_3}(z,t/d)+dt^2$ is also $G$-invariant and we now
glue to get a smooth $G$-invariant metric with scalar curvature
bigger than $s_g-\epsilon_2$ on $M$.

An important fact about the bending of $\gamma$ is that if we can
take $r_1'$ and $r_2$ further small, we only need to shrink the
remaining part of $\gamma$ homothetically. Let $\{(t,f(t))\}$ be
the graph of $\gamma$ in step 3 and $\tau_1$ be $f^{-1}(r_2)$. For
$0 < \mu \leq 1$, let's take $\mu r_1'$ and $\mu r_2$ instead of
$r_1'$ and $r_2$ respectively, and let $\tau_{\mu}$ be the
$t$-coordinate corresponding to $\mu r_2$. Then we shrink the step
3 part of $\gamma$ homothetically by $\mu$ and concatenate it to
$(\tau_{\mu}, \mu r_2)$. Indeed the equation of this portion of
the curve is given by $(t,\mu
f(\frac{t-\tau_{\mu}+\mu\tau_1}{\mu}))$. Moreover, noting that the
geodesic curvature $k$ is dilated by $\frac{1}{\mu}$ without
changing $\theta$, the scalar curvature at $(t,\mu y)$ satisfies
\begin{eqnarray*}
s(t,\mu y) &\geq& s_g(t,\mu y)+
\frac{(q-1)(q-2)}{2}\frac{\sin^2\theta}{(\mu |y|)^2}
-3(q-1)\frac{k\sin\theta}{\mu |y|}\\ &\geq& s_g(t,\mu y),
\end{eqnarray*}
where we used (\ref{original}) in the second inequality. We denote
the curve with $\mu r_1'$ and $\mu r_2$ instead of $r_1'$ and
$r_2$ by $\gamma_{\mu}$.

We also claim that the metric on the homotopy region $W \times
S^{q-1} \times [0,d]$ can be accordingly shrunk to $H_{\mu
r_3}(z,t/d)+\mu^2 dt^2$ still having positive scalar curvature for
any $\mu \in (0,1]$, once $r_2$ and hence $r_3$ was chosen
sufficiently small.

\begin{lem}
The scalar curvature of the manifold $W \times S^{q-1} \times
[0,d]$ with the metric $H_{\mu r_3}(z,t/d)+\mu^2 dt^2$ is bounded
below by $\frac{C}{(\mu r_3)^2}+\frac{C'}{\mu^2}$ for any $\mu\in
(0,1]$, and any sufficiently small $r_3>0$, where $C>0$ is given
in lemma \ref{zz} and $C'$ is a constant.
\end{lem}
\begin{proof}
The proof continues from the above lemma. Using the estimates
(\ref{mat1}) and (\ref{mat2}), $H_{\mu r_3}(z,t/d)+\mu^2 dt^2$ is
given by
$$
\left(
\begin{tabular}{c|c}
&\\
$ $&$ $\\
&\\
$\ \ \ H_{\mu
r_3}(z,t/d)\ \ \ $&$0$\\
&\\
$ $&$ $\\
&\\
\cline{1-2}&\\
$0$ & $\mu^2$\\&\\
\end{tabular}
 \right)
 =
\left(
\begin{tabular}{c|c}
&\\
\begin{tabular}{c|c}
&\\
$O(1)$&$O((\mu r_3)^2)$\\ &\\
\cline{1-2}&\\
$O((\mu r_3)^2)$ & $O((\mu r_3)^2)$\\&\\
\end{tabular}
&$0$\\ &\\
\cline{1-2}&\\
$0$ & $\mu^2$\\&\\
\end{tabular}
 \right),
 $$ and its inverse is given by
$$
\left(
\begin{tabular}{c|c}
&\\
$ $&$ $\\
&\\
$(H_{\mu
r_3}(z,t/d))^{-1}$&$0$\\
&\\
$ $&$ $\\
&\\
\cline{1-2}&\\
$0$ & $\frac{1}{\mu^2}$\\&\\
\end{tabular}
 \right)
 =
\left(
\begin{tabular}{c|c}
&\\
\begin{tabular}{c|c}
&\\
$O(1)$&$O(1)$\\ &\\
\cline{1-2}&\\
$O(1)$ & $O(\frac{1}{(\mu r_3)^2})$\\&\\
\end{tabular}
&$0$\\ &\\
\cline{1-2}&\\
$0$ & $\frac{1}{\mu^2}$\\&\\
\end{tabular}
 \right).
 $$
The same estimates also hold for their derivatives. We let
$\Gamma^{\mu}$ and $R^{\mu}$ be the Christoffel symbol and the
Riemann curvature tensor of $H_{\mu r_3}(z,t/d)+\mu^2 dt^2$
respectively. As before $A,B,C,\cdots$ run from $1$ to $n-1$, and
$N$ denotes the index of the last coordinate function $t$. The
direct computations show that
$$(\Gamma^\mu)_{NN}^{N}=(\Gamma^\mu)_{AN}^{N}=0=
\partial(\Gamma^\mu)_{NN}^{N}=\partial(\Gamma^\mu)_{AN}^{N},$$
$$(\Gamma^\mu)_{AB}^{N}=\frac{1}{\mu^2}O(1)
=\partial(\Gamma^\mu)_{AB}^{N},\ \ \ \
(\Gamma^\mu)_{AN}^{C}=O(1)=\partial(\Gamma^\mu)_{AN}^{C},$$ and
$$(R^\mu)_{NBC}^{N}=\frac{1}{\mu^2}O(1).$$

Let $X_t$ be the hypersurface $W\times S^{q-1} \times \{t\}$. Then
the second fundamental form of $X_t$ is given by
$(\Gamma^\mu)_{AB}^{N}=\frac{1}{\mu^2}O(1)$, and hence its norm is
of the form $\frac{1}{\mu}O(1)$. Denote the scalar curvature of
the hypersurface $X_t$ with the induced metric by $s_{X_t}$. It
follows from the Gauss curvature equation and the above lemma that
the scalar curvature is given by
\begin{eqnarray*}
s_{X_t}+\frac{1}{\mu^2}O(1) +2\sum_{B=1}^{n}(R^\mu)_{NBB}^{N}
&=&s_{X_t}+\frac{1}{\mu^2}O(1)\\ &\geq& \frac{C}{(\mu
r_3)^2}+\frac{1}{\mu^2}O(1).
\end{eqnarray*}
\end{proof}

Therefore the scalar curvature of $H_{\mu r_3}(z,t/d)+\mu^2 dt^2$
is positive for sufficiently small $r_3>0$. From now on we assume
that $r_2$ was taken small enough to ensure this, and the
Riemannian $G$-manifold obtained by $\gamma_\mu$ and $H_{\mu
r_3}(z,t/d)+\mu^2 dt^2$ is denoted by $(M_\mu, \tilde{g}_\mu)$.

We define three Riemannian manifolds with boundary
$(S_{\delta,\varepsilon},\tilde{g}_{\delta,\varepsilon}) \subset
(T_{\delta,\varepsilon},\tilde{g}_{\delta,\varepsilon}) \subset
(N_{\delta\varepsilon},\tilde{g}_{\delta\varepsilon})$ by
$$S_{\delta,\varepsilon}\equiv M_{\delta\varepsilon}-(M_0-\{r\leq
\delta\varepsilon r_1'\}),$$
$$T_{\delta,\varepsilon}\equiv M_{\delta\varepsilon}-(M_0-\{r\leq
\varepsilon r_1\}),$$ and $$N_{\delta\varepsilon}\equiv
M_{\delta\varepsilon}-(M_0-\{r\geq r_0\})$$ with the induced
metric. (In fact,
$(S_{\delta,\varepsilon},\tilde{g}_{\delta,\varepsilon})$ depends
only on $\delta\varepsilon$.) To investigate the relation between
$T_{\delta,1}$ and $T_{\delta,\varepsilon}$, let $x$ be any point
in $W$ and define a $q$-dimensional Riemanian submanifold
$(T_{\delta,\varepsilon,x},\tilde{g}_{\delta,\varepsilon,x})\subset
(T_{\delta,\varepsilon},\tilde{g}_{\delta,\varepsilon})$ by
$T_{\delta,\varepsilon,x}\equiv T_{\delta,\varepsilon}\cap (\{x\}
\times S^{q-1}\times [0,d])$ with the induced metric. Taking into
account that $g$ is $C^0$-near to the product metric on $N(r_0)$,
i.e. $g=g^W+ g_{E}+O(r_0)$, where $g_{E}$ is the Euclidean metric
on $\Bbb R^q$, we have
$$\tilde{g}_{\delta,\varepsilon}=g^W+\tilde{g}_{\delta,\varepsilon,x}+O(\varepsilon
r_1)$$ on $T_{\delta,\varepsilon}$. The obvious shrinking map from
$\gamma_{\delta}$ for $r\leq r_1$ onto
$\gamma_{\delta\varepsilon}$ for $r\leq \varepsilon r_1$ and the
identity map in the homotopy region induces a diffeomorphism
$\Phi_{\delta,\varepsilon}$ from $T_{\delta,1}$ to
$T_{\delta,\varepsilon}$, which gives
$\Phi^*(\tilde{g}_{\delta,\varepsilon,x})=\varepsilon^2
\tilde{g}_{\delta,1,x}$. Thus we have on $T_{\delta,\varepsilon}$,
\begin{eqnarray} \label{vol}
\Phi^*(dV_{\tilde{g}_{\delta,\varepsilon}})&=&\Phi^*((1+O(\varepsilon
r_1))dV_{g^W} dV_{\tilde{g}_{\delta,\varepsilon,x}})=\varepsilon^q
(1+O(\varepsilon r_1))dV_{g^W} dV_{\tilde{g}_{\delta,1,x}}\nonumber \\
&\lessgtr&\varepsilon^q (1\pm C_1r_1)dV_{\tilde{g}_{\delta,1}},
\end{eqnarray}
where $C_1 >0$ is a constant. From now on $C_i$'s will denote some
positive constants. Let $\langle \cdot, \cdot
\rangle_{\tilde{g}_{\delta,\varepsilon}}$ and $\langle \cdot,
\cdot \rangle_{\tilde{g}_{\delta,\varepsilon,x}}$ denote the inner
product on
$(T_{\delta,\varepsilon},\tilde{g}_{\delta,\varepsilon})$ and
$(T_{\delta,\varepsilon,x},\tilde{g}_{\delta,\varepsilon,x})$
respectively. Then we also have on $T_{\delta,\varepsilon}$,
\begin{eqnarray} \label{form}
\Phi^*\langle\omega,\omega\rangle_{\tilde{g}_{\delta,\varepsilon}}&=&\Phi^*
\langle\omega,\omega\rangle_{\tilde{g}_{\delta,\varepsilon,x}}=\frac{1}{\varepsilon^2}
\langle \omega, \omega \rangle_{\tilde{g}_{\delta,1,x}}\nonumber \\
&=& \frac{1}{\varepsilon^2}
\langle\omega,\omega\rangle_{\tilde{g}_{\delta,1}}
\end{eqnarray}
for any $1$-form $\omega$ belonging to $T^*(S^{q-1} \times [0,d])$
in $T^*(W \times S^{q-1} \times [0,d])$. It's important that $C_1$
is a uniform constant independent of any choices we made such as
$\theta_0,r_2,\delta,$ and etc, as long as $r_0$ is sufficiently
small, which we always assume. From now on we will omit $\Phi^*$
for convenience. Also note that for any choice of $r_0$ and
$\theta_0$, the length of the step 3 part of
$\gamma_{\delta\varepsilon}$ and the volume of the homotopy region
can be made arbitrarily small by taking $r_2$ much smaller, which
we always assume from now on. This means that there exist
constants $C_3,C_4,C_5>0$ such that
$$\textrm{vol}_{\tilde{g}_{\delta\varepsilon}}(N_{\delta\varepsilon})\leq
C_3r_0^q,\ \ \ \
\textrm{vol}_{\tilde{g}_{\delta,\varepsilon}}(T_{\delta,\varepsilon})\geq
C_4(\varepsilon r_1)^q,$$ and
$$\textrm{vol}_{\tilde{g}_{\delta,\varepsilon}}(S_{\delta,\varepsilon})\leq
C_5(\delta\varepsilon r_1')^{q},$$ where $C_i$'s are also uniform
constants when $r_0,\theta_0$ and $r_2$ are chosen small by the
above way. As the last preparation, we have
\begin{lem}
There is a constant $\hat{C}>0$ independent of $\delta\in (0,1]$
satisfying the Sobolev inequality
\begin{equation}\label{sobo}
(\frac{\int_{S_{\delta,1}} \varphi^p\
dV_{\tilde{g}_{\delta,1}}}{\textrm{vol}_{\tilde{g}_{\delta,1}}(S_{\delta,1})}
)^{\frac{1}{p}}\leq \hat{C}((\frac{\int_{S_{\delta,1}} \varphi^2\
dV_{\tilde{g}_{\delta,1}}}{\textrm{vol}_{\tilde{g}_{\delta,1}}(S_{\delta,1})}
)^{\frac{1}{2}}+(\frac{\int_{S_{\delta,1}}
|d\varphi|^2_{\tilde{g}_{\delta,1}}
dV_{\tilde{g}_{\delta,1}}}{\textrm{vol}_{\tilde{g}_{\delta,1}}(S_{\delta,1})})^{\frac{1}{2}})
\end{equation}
for any $\varphi\in L^2_1(S_{\delta,1})$.
\end{lem}
\begin{proof}
For a fixed $\theta_0,r_1',$ and $r_2$, get
$(S_{1,1},\tilde{g}_{1,1})$ and choose a $\hat{C}$ satisfying the
above inequality. In the same way as above, consider a
diffeomorphism $\Psi$ from $S_{1,1}$ onto $S_{\delta,1}$ such that
$$\Psi^*(dV_{\tilde{g}_{\delta,1}})\lessgtr\delta^q
(1\pm C_6r_1')dV_{\tilde{g}_{1,1}},$$
$$\Psi^*\langle\omega,\omega\rangle_{\tilde{g}_{\delta,1}}=\frac{1}{\delta^2}
\langle\omega,\omega\rangle_{\tilde{g}_{1,1}},$$ and
$$\Psi^*\langle\sigma,\sigma\rangle_{\tilde{g}_{\delta,1}}\lessgtr
(1\pm C_8r_1')\langle\sigma,\sigma\rangle_{\tilde{g}_{1,1}}$$ for
any $1$-forms $\omega$ and $\sigma$ belonging to $T^*(S^{q-1}
\times [0,d])$ and $T^*W$ in $T^*(W \times S^{q-1} \times [0,d])$
respectively. Then the result follows immediately.
\end{proof}

\noindent Although it is not necessary for our further discussion,
we remark that

\begin{rmk}
In fact $\hat{C}$ may depend only on $\theta_0,r_1',$ and $r_2.$
Notice that $\hat{C}$ is a continuous function of the metric in
$C^0$-norm. Since the ambiguity of the step 3 construction of
$\gamma$ can be made very small, any possible
$(S_{1,1},\tilde{g}_{1,1})$ is $C^0$-close, once $\theta_0,r_1',$
$r_2$ are determined. As a final note, actually we will not need
the $\delta$-independence of $\hat{C}$, because we will use
$\hat{C}$ for a fixed $\delta$.
\end{rmk}

Now let's get down to estimating the $G$-Yamabe constant of
$(M_{\delta\varepsilon},[\tilde{g}_{\delta\varepsilon}]_G)$. Let
$\varphi_{\delta\varepsilon}$ be a $G$-Yamabe minimizer satisfying
$\int_{M_{\delta\varepsilon}} \varphi^p_{\delta\varepsilon}\
dV_{\tilde{g}_{\delta\varepsilon}}=1.$ We have two cases, either
$$\int_{S_{\delta,\varepsilon}}\varphi^p_{\delta\varepsilon}\
dV_{\tilde{g}_{\delta,\varepsilon}} \leq \frac{2^{p+1}\hat{C}^p\
\textrm{vol}_{\tilde{g}_{\delta,1}}(S_{\delta,1})}
{\textrm{vol}_{\tilde{g}_{\delta,1}}(T_{\delta,1}-S_{\delta,1})}
\int_{M_{\delta\varepsilon}-S_{\delta,\varepsilon}}\varphi^p_{\delta\varepsilon}\
dV_{\tilde{g}_{\delta\varepsilon}}$$ or not.

Assume the first case. Let $\eta_{\delta\varepsilon}(r)$ be
defined by $\eta(\frac{r}{\delta\varepsilon})$. On the support of
$\eta_{\delta\varepsilon}$, $\tilde{g}_{\delta\varepsilon}$ is
very close to $g$ when $\theta_0$ is very small. To compare these
two metrics on this region, let $i: M_0- N(\delta\varepsilon r_3)
\rightarrow M_{\delta\varepsilon}$ be the obvious inclusion map.
Then $i$ is isometric on the outside of $N(r_0)$. On
$N(r_0)-N(\delta\varepsilon r_3)$, $i$ is isometric in the
direction orthogonal to the radial direction, and
$\frac{\partial}{\partial r}$ gets dilated by
$\frac{1}{\sqrt{1-\sin^2 \theta}}$. In particular on the support
of $\eta_{\delta\varepsilon}$,
$$dV_{\tilde{g}_{\delta\varepsilon}} \geq dV_{g}\geq
\sqrt{1-\sin^2 \theta_0}\ dV_{\tilde{g}_{\delta\varepsilon}},$$
and
$$|\omega|_{\tilde{g}_{\delta\varepsilon}} \leq |\omega|_{g}\leq
\frac{1}{\sqrt{1-\sin^2 \theta_0}}\
|\omega|_{\tilde{g}_{\delta\varepsilon}}$$ for any $1$-form
$\omega$. This gives us that
$$\int_{M_{\delta\varepsilon}} \varphi^p_{\delta\varepsilon}\
dV_{\tilde{g}_{\delta\varepsilon}} \geq \int_{M_0}
(\eta_{\delta\varepsilon}\varphi_{\delta\varepsilon})^p\ dV_{g},$$
and
\begin{eqnarray*}
\int_{M_0}
(\eta_{\delta\varepsilon}\varphi_{\delta\varepsilon})^p\
dV_{g}&\geq& \sqrt{1-\sin^2 \theta_0} \int_{M_{\delta\varepsilon}}
(\eta_{\delta\varepsilon}\varphi_{\delta\varepsilon})^p\
dV_{\tilde{g}_{\delta\varepsilon}} \\ &\geq& \sqrt{1-\sin^2
\theta_0}(\int_{M_{\delta\varepsilon}}
\varphi^p_{\delta\varepsilon}\
dV_{\tilde{g}_{\delta\varepsilon}}-\int_{S_{\delta,\varepsilon}}
\varphi^p_{\delta\varepsilon}\
dV_{\tilde{g}_{\delta,\varepsilon}})\\
&\geq& \sqrt{1-\sin^2 \theta_0}(\int_{M_{\delta\varepsilon}}
\varphi^p_{\delta\varepsilon}\ dV_{\tilde{g}_{\delta\varepsilon}}\\
& & -\frac{2^{p+1}\hat{C}^p\
\textrm{vol}_{\tilde{g}_{\delta,1}}(S_{\delta,1})}
{\textrm{vol}_{\tilde{g}_{\delta,1}}(T_{\delta,1}-S_{\delta,1})}
\int_{M_{\delta\varepsilon}-S_{\delta,\varepsilon}}\varphi^p_{\delta\varepsilon}\
dV_{\tilde{g}_{\delta\varepsilon}})\\
&\geq& \sqrt{1-\sin^2
\theta_0}((1-\frac{2^{p+1}\hat{C}^pC_5(\delta r_1')^q}{C_4
r_1^q-C_5 (\delta r_1')^q})\int_{M_{\delta\varepsilon}}
\varphi^p_{\delta\varepsilon}\ dV_{\tilde{g}_{\delta\varepsilon}})\\
&=&\sqrt{1-\sin^2 \theta_0}
(1-C_9\delta^q)\int_{M_{\delta\varepsilon}}
\varphi^p_{\delta\varepsilon}\
dV_{\tilde{g}_{\delta\varepsilon}}).
\end{eqnarray*}
Using the fact that $s_{\tilde{g}_{\delta\varepsilon}} \geq
s_g+\frac{(q-1)(q-2)}{2}\frac{\sin^2 \theta_0}{r^2}\geq
s_g+|d\eta_{\delta\varepsilon}|_{g}^2$ on the support of
$d\eta_{\delta\varepsilon}$, and
$s_{\tilde{g}_{\delta\varepsilon}}$ is bounded below by $(\min
s_g)-\epsilon_2$, we get
\begin{align*}
&\int_{M_0}(|d(\eta_{\delta\varepsilon}\varphi_{\delta\varepsilon})|_{g}^2
+s_{g}(\eta_{\delta\varepsilon}\varphi_{\delta\varepsilon})^2)\ dV_{g}\\
&=\int_{\{r\geq \delta\varepsilon
r_2\}}(\eta_{\delta\varepsilon}^2|d\varphi_{\delta\varepsilon}|_{g}^2
+|d\eta_{\delta\varepsilon}|_{g}^2\varphi_{\delta\varepsilon}^2 +
s_{g}(\eta_{\delta\varepsilon}\varphi_{\delta\varepsilon}^2))\ dV_{g}\\
&\leq\int_{M_{\delta\varepsilon}}\frac{1}{1-\sin^2
\theta_0}(|d\varphi_{\delta\varepsilon}|_{\tilde{g}_{\delta\varepsilon}}^2
+s_{\tilde{g}_{\delta\varepsilon}} \varphi_{\delta\varepsilon}^2)\
dV_{\tilde{g}_{\delta\varepsilon}}\\ & +
C_{10}\int_{\{\delta\varepsilon r_2\leq r\leq \delta\varepsilon
r_1'\}} \varphi^2_{\delta\varepsilon}\
dV_{\tilde{g}_{\delta\varepsilon}}+ C_{11}(\sin^2
\theta_0+\epsilon_2)\int_{M_{\delta\varepsilon}}
\varphi^2_{\delta\varepsilon}\ dV_{\tilde{g}_{\delta\varepsilon}},
\end{align*}
where $C_{10}$ and $C_{11}$ are constants depending only on $\min
s_g$. By using the H\"{o}lder inequality the second term is
bounded above by
$$C_{10}
(\textrm{vol}(S_{\delta\varepsilon})_{\tilde{g}_{\delta\varepsilon}})^{\frac{2}{n}}
(\int_{S_{\delta\varepsilon}} \varphi^p_{\delta\varepsilon}\
dV_{\tilde{g}_{\delta\varepsilon}})^{\frac{2}{p}}\leq C_{10}(C_5
(\delta\varepsilon
r_1')^q)^{\frac{2}{n}}(\int_{M_{\delta\varepsilon}}
\varphi^p_{\delta\varepsilon}\
dV_{\tilde{g}_{\delta\varepsilon}})^{\frac{2}{p}},$$ and the third
term is bounded above by
$$C_{11} (\sin^2 \theta_0+\epsilon_2)
(\textrm{vol}(M_{\delta\varepsilon})_{\tilde{g}_{\delta\varepsilon}})^{\frac{2}{n}}
(\int_{M_{\delta\varepsilon}} \varphi^p_{\delta\varepsilon}\
dV_{\tilde{g}_{\delta\varepsilon}})^{\frac{2}{p}}\leq
C_{12}(\sin^2 \theta_0+\epsilon_2)(\int_{M_{\delta\varepsilon}}
\varphi^p_{\delta\varepsilon}\
dV_{\tilde{g}_{\delta\varepsilon}})^{\frac{2}{p}},$$ where we used
the fact that
\begin{eqnarray*}
\textrm{vol}_{\tilde{g}_{\delta\varepsilon}}(M_{\delta\varepsilon})&=&
\textrm{vol}_{g}(M_0-N(r_0))+\textrm{vol}_{\tilde{g}_{\delta\varepsilon}}
(N_{\delta\varepsilon})\\ &\leq&
\textrm{vol}_{g}(M_0-N(r_0))+C_{3}r_0^q.
\end{eqnarray*}
Thus
$$ Q_g(\eta_{\delta\varepsilon}\varphi_{\delta\varepsilon})\geq Y(M_0,[g]_G)\geq
Y_G(M_0)-\epsilon_1$$ is bounded above by
$$ (1-\sin^2\theta_0)^{-1}Y(M_{\delta\varepsilon},[\tilde{g}_{\delta\varepsilon}]_G)
+C_{10}(C_5 (\delta\varepsilon
r_1')^q)^{\frac{2}{n}}+C_{12}(\sin^2 \theta_0+\epsilon_2)$$ or
$$\frac{(1-\sin^2\theta_0)^{-1}Y(M_{\delta\varepsilon},[\tilde{g}_{\delta\varepsilon}]_G)
+C_{10}(C_5 (\delta\varepsilon
r_1')^q)^{\frac{2}{n}}+C_{12}(\sin^2 \theta_0+\epsilon_2)}
{(1-\sin^2 \theta_0)^{\frac{1}{p}}(1-C_9\delta^q)^{\frac{2}{p}}}$$
for any $\delta$ and $\varepsilon$. Recall that
$C_{10}C_5^{\frac{2}{n}}$ and $C_{12}$ are uniform constants
independent of any choices and $C_9$ is independent of $\delta$
and $\varepsilon$. Taking first $\theta_0$ and then $\delta$
arbitrarily small, we have
$$Y_G(M)+C_{12}\epsilon_2\geq Y_G(M_0)-\epsilon_1.$$ Since $\epsilon_1,\epsilon_2>0$ are arbitrary,
it follows that
$$Y_G(M)\geq Y_G(M_0).$$

In the second case, we want to derive a contradiction when
$\varepsilon >0$ gets sufficiently small for any fixed $\delta>0$.

\begin{lem}
Suppose that $(X,h)$ is a compact Riemannian manifold with smooth
boundary and that $f \in L^2(X)$ satisfying $\int_X f\ dV_h=0$.
Then there exists a function $\xi \in L^2_2(X)$ unique up to the
addition of constant such that $\Delta \xi = f$ and in addition
$\vec{n}\cdot \nabla\xi$ vanishes at the boundary, where $\vec{n}$
is the unit outward normal to the boundary.
\end{lem}
\begin{proof}
See \cite{Hor}.
\end{proof}

Consider a step function $f_\delta$ on $T_{\delta,1}$ defined by
$$f_\delta=\left\{ \begin{array}{cl}
\textrm{vol}_{\tilde{g}_{\delta,1}}(S_{\delta,1})^{-1}&\mbox{on } S_{\delta,1}\\
(\textrm{vol}_{\tilde{g}_{\delta,1}}(S_{\delta,1})-
\textrm{vol}_{\tilde{g}_{\delta,1}}(T_{\delta,1}))^{-1}& \mbox{on
} T_{\delta,1}-S_{\delta,1}.
\end{array}\right.$$
Then $\int_{T_{\delta,1}} f_\delta\ dV_{\tilde{g}_{\delta,1}}=0$,
so by the above lemma, there exists a function $\xi_\delta \in
L^2_2(T_{\delta,1})$ satisfying $\Delta\xi_\delta=f_\delta$, and
that $\nabla \xi_\delta$ vanishes normal to the boundary. For any
$\varphi \in L^2_1(T_{\delta,1})$ the integration by parts yields
\begin{align*}
\frac{1}{\textrm{vol}_{\tilde{g}_{\delta,1}}(S_{\delta,1})}&
\int_{S_{\delta,1}} \varphi\ dV_{\tilde{g}_{\delta,1}} -
\frac{1}{\textrm{vol}_{\tilde{g}_{\delta,1}}(T_{\delta,1}-S_{\delta,1})}
\int_{T_{\delta,1}-S_{\delta,1}}\varphi\ dV_{\tilde{g}_{\delta,1}} \\
&= \int_{T_{\delta,1}}\varphi \Delta\xi_\delta\
dV_{\tilde{g}_{\delta,1}}=\int_{T_{\delta,1}}\langle d \varphi,d
\xi_\delta\rangle_{\tilde{g}_{\delta,1}}
dV_{\tilde{g}_{\delta,1}},
\end{align*}
and hence
\begin{align}\label{1}
\frac{1}{\textrm{vol}_{\tilde{g}_{\delta,1}}(S_{\delta,1})}& |
\int_{S_{\delta,1}} \varphi\ dV_{\tilde{g}_{\delta,1}}| -
\frac{1}{\textrm{vol}_{\tilde{g}_{\delta,1}}(T_{\delta,1}-S_{\delta,1})}
|\int_{T_{\delta,1}-S_{\delta,1}}\varphi\
dV_{\tilde{g}_{\delta,1}}| \nonumber
\\ &\leq (\int_{T_{\delta,1}} |d
\xi_\delta|^2_{\tilde{g}_{\delta,1}}
dV_{\tilde{g}_{\delta,1}})^{\frac{1}{2}} (\int_{T_{\delta,1}} |d
\varphi|^2_{\tilde{g}_{\delta,1}}
dV_{\tilde{g}_{\delta,1}})^{\frac{1}{2}}
\end{align}
by the H\"{o}lder inequality. Since $S_{\delta,1}$ is connected,
the constants are the only eigenvectors of $\Delta$ on
$S_{\delta,1}$ with eigenvalue $0$ and derivative vanishing normal
to the boundary. By the discreteness of the spectrum of $\Delta$
on $S_{\delta,1}$ with these boundary conditions, we have
\begin{equation}\label{2}
(\frac{\int_{S_{\delta,1}} \varphi^2\
dV_{\tilde{g}_{\delta,1}}}{\textrm{vol}_{\tilde{g}_{\delta,1}}(S_{\delta,1})}
)^{\frac{1}{2}}\leq
\frac{1}{\textrm{vol}_{\tilde{g}_{\delta,1}}(S_{\delta,1})}|\int_{S_{\delta,1}}
\varphi\ dV_{\tilde{g}_{\delta,1}}|+(\frac{\int_{S_{\delta,1}}
|d\varphi|^2_{\tilde{g}_{\delta,1}}
dV_{\tilde{g}_{\delta,1}}}{C_{13}\textrm{vol}_{\tilde{g}_{\delta,1}}(S_{\delta,1})}
)^{\frac{1}{2}}.
\end{equation}
Also, by the Sobolev inequality (\ref{sobo}),
\begin{equation}\label{3}
\frac{1}{\hat{C}}(\frac{\int_{S_{\delta,1}} \varphi^p\
dV_{\tilde{g}_{\delta,1}}}{\textrm{vol}_{\tilde{g}_{\delta,1}}(S_{\delta,1})}
)^{\frac{1}{p}}-(\frac{\int_{S_{\delta,1}}
|d\varphi|^2_{\tilde{g}_{\delta,1}}
dV_{\tilde{g}_{\delta,1}}}{\textrm{vol}_{\tilde{g}_{\delta,1}}
(S_{\delta,1})})^{\frac{1}{2}}\leq(\frac{\int_{S_{\delta,1}}
\varphi^2\
dV_{\tilde{g}_{\delta,1}}}{\textrm{vol}_{\tilde{g}_{\delta,1}}(S_{\delta,1})}
)^{\frac{1}{2}}.
\end{equation}
On the other hand, the H\"{o}lder inequality gives
\begin{equation}\label{4}
\frac{1}{\textrm{vol}_{\tilde{g}_{\delta,1}}(T_{\delta,1}-S_{\delta,1})}
|\int_{T_{\delta,1}-S_{\delta,1}}\varphi\
dV_{\tilde{g}_{\delta,1}}|
\leq(\frac{\int_{T_{\delta,1}-S_{\delta,1}} \varphi^p\
dV_{\tilde{g}_{\delta,1}}}{\textrm{vol}_{\tilde{g}_{\delta,1}}
(T_{\delta,1})-\textrm{vol}_{\tilde{g}_{\delta,1}}(S_{\delta,1})}
)^{\frac{1}{p}}.
\end{equation}
Adding together (\ref{1}), (\ref{2}), (\ref{3}), and (\ref{4})
yields
\begin{align*}
\frac{1}{\hat{C}}(&\frac{\int_{S_{\delta,1}} \varphi^p\
dV_{\tilde{g}_{\delta,1}}}{\textrm{vol}_{\tilde{g}_{\delta,1}}(S_{\delta,1})}
)^{\frac{1}{p}} -(\frac{\int_{T_{\delta,1}-S_{\delta,1}}
\varphi^p\
dV_{\tilde{g}_{\delta,1}}}{\textrm{vol}_{\tilde{g}_{\delta,1}}
(T_{\delta,1})-\textrm{vol}_{\tilde{g}_{\delta,1}}(S_{\delta,1})}
)^{\frac{1}{p}}\\ &{\ \ \ \ \ \ \ \ \ \ } \leq
C_{14}(\int_{T_{\delta,1}} |d \varphi|^2_{\tilde{g}_{\delta,1}}
dV_{\tilde{g}_{\delta,1}})^{\frac{1}{2}}.
\end{align*}

Now if $\varphi$ is $G$-invariant, then
$\frac{\partial\varphi}{\partial x_i}=0$ for any $i=1,\cdots,n-q$,
because the $G$-action on $W$ is locally transitive. Then using
(\ref{vol}) and (\ref{form}), we get
\begin{align*}
\frac{1}{\hat{C}}(\frac{1-C_{15} r_1}{\varepsilon^q}
&\frac{\int_{S_{\delta,\varepsilon}}\varphi^p\
dV_{\tilde{g}_{\delta,\varepsilon}}}{\textrm{vol}_{\tilde{g}_{\delta,1}}
(S_{\delta,1})})^{\frac{1}{p}}-(\frac{1+C_{16} r_1}{\varepsilon^q}
\frac{\int_{T_{\delta,\varepsilon}-S_{\delta,\varepsilon}}\varphi^p\
dV_{\tilde{g}_{\delta,\varepsilon}}}{\textrm{vol}_{\tilde{g}_{\delta,1}}
(T_{\delta,1}-S_{\delta,1})})^{\frac{1}{p}}\\
&\leq
C_{14}(\frac{1+C_{17}r_1}{\varepsilon^{q-2}}\int_{T_{\delta,\varepsilon}}
|d\varphi|^2_{\tilde{g}_{\delta,\varepsilon}}
dV_{\tilde{g}_{\delta,\varepsilon}})^{\frac{1}{2}}.
\end{align*}
Since $C_{15},C_{16},$ and $C_{17}$ are uniform constants, we get
for sufficiently small $r_1>0$,
\begin{align*}
(\int_{S_{\delta,\varepsilon}} \varphi^p\
dV_{\tilde{g}_{\delta,\varepsilon}})^{\frac{1}{p}}&-(\frac{2\
\hat{C}^p\
\textrm{vol}_{\tilde{g}_{\delta,1}}(S_{\delta,1})}{\textrm{vol}_{\tilde{g}_{\delta,1}}
(T_{\delta,1}-S_{\delta,1})}
\int_{T_{\delta,\varepsilon}-S_{\delta,\varepsilon}}
\varphi^p\ dV_{\tilde{g}_{\delta,\varepsilon}})^{\frac{1}{p}}\\
&\leq C_{18}\ \varepsilon^{1-\frac{q}{n}}
(\int_{T_{\delta,\varepsilon}} |d
\varphi|^2_{\tilde{g}_{\delta,\varepsilon}}
dV_{\tilde{g}_{\delta,\varepsilon}})^{\frac{1}{2}}.
\end{align*}

Under the assumption that
$$\int_{S_{\delta,\varepsilon}}\varphi_{\delta\varepsilon}^p\
dV_{\tilde{g}_{\delta,\varepsilon}} > \frac{2^{p+1} \hat{C}^p\
\textrm{vol}_{\tilde{g}_{\delta,1}}(S_{\delta,1})}{\textrm{vol}_{\tilde{g}_{\delta,1}}
(T_{\delta,1}-S_{\delta,1})}
\int_{M_{\delta\varepsilon}-S_{\delta,\varepsilon}}\varphi_{\delta\varepsilon}^p\
dV_{\tilde{g}_{\delta\varepsilon}},$$ we have
$$\frac{1}{2}(\int_{S_{\delta,\varepsilon}} \varphi_{\delta\varepsilon}^p\
dV_{\tilde{g}_{\delta,\varepsilon}})^{\frac{1}{p}} \leq C_{18}\
\varepsilon^{1-\frac{q}{n}} (\int_{T_{\delta,\varepsilon}} |d
\varphi_{\delta\varepsilon}|^2_{\tilde{g}_{\delta,\varepsilon}}
dV_{\tilde{g}_{\delta,\varepsilon}})^{\frac{1}{2}},$$ and
\begin{eqnarray*}
(\int_{M_{\delta\varepsilon}} \varphi_{\delta\varepsilon}^p\
dV_{\tilde{g}_{\delta\varepsilon}})^{\frac{1}{p}}&=&(\int_{S_{\delta,\varepsilon}}
\varphi_{\delta\varepsilon}^p\
dV_{\tilde{g}_{\delta,\varepsilon}}+\int_{M_{\delta\varepsilon}-S_{\delta,\varepsilon}}
\varphi_{\delta\varepsilon}^p\ dV_{\tilde{g}_{\delta\varepsilon}}
)^{\frac{1}{p}}
\\ &\leq& (\int_{S_{\delta,\varepsilon}} \varphi_{\delta\varepsilon}^p\
dV_{\tilde{g}_{\delta,\varepsilon}})^{\frac{1}{p}}+
(\int_{M_{\delta\varepsilon}-S_{\delta,\varepsilon}}
\varphi_{\delta\varepsilon}^p\ dV_{\tilde{g}_{\delta\varepsilon}})^{\frac{1}{p}} \\
&\leq& (\int_{S_{\delta,\varepsilon}}
\varphi_{\delta\varepsilon}^p\
dV_{\tilde{g}_{\delta,\varepsilon}})^{\frac{1}{p}}(1+
(\frac{\textrm{vol}_{\tilde{g}_{\delta,1}}
(T_{\delta,1}-S_{\delta,1})}{2^{p+1} \hat{C}^p\
\textrm{vol}_{\tilde{g}_{\delta,1}}(S_{\delta,1})})^{\frac{1}{p}}),
\end{eqnarray*}
which yield $$C_{19}\ \varepsilon^{-2+\frac{2q}{n}}\leq
\frac{\int_{T_{\delta,\varepsilon}}
|d\varphi_{\delta\varepsilon}|^2_{\tilde{g}_{\delta,\varepsilon}}
dV_{\tilde{g}_{\delta,\varepsilon}}}{(\int_{M_{\delta\varepsilon}}
\varphi_{\delta\varepsilon}^p\
dV_{\tilde{g}_{\delta\varepsilon}})^{\frac{2}{p}}}.$$

On the other hand, using the fact that
$s_{\tilde{g}_{\delta\varepsilon}}$ is bounded below and
$\textrm{vol}_{\tilde{g}_{\delta\varepsilon}}(M_{\delta\varepsilon})$
is bounded above for any $\varepsilon\in (0,1]$, a simple
application of the H\"{o}lder inequality gives
\begin{eqnarray*}
\frac{\int_{M_{\delta\varepsilon}}
s_{\tilde{g}_{\delta\varepsilon}} \varphi_{\delta\varepsilon}^2\
dV_{\tilde{g}_{\delta\varepsilon}}}{(\int_{M_{\delta\varepsilon}}
\varphi_{\delta\varepsilon}^p\
dV_{\tilde{g}_{\delta\varepsilon}})^{\frac{2}{p}}}&\geq&
\min(0,(\min
s_g)-\epsilon_2)(\textrm{vol}_{\tilde{g}_{\delta\varepsilon}}
(M_{\delta\varepsilon}))^{1-\frac{2}{p}}\\ &\geq& -C_{20}
\end{eqnarray*}
for any $\varepsilon$. Now by letting $\varepsilon \rightarrow 0$,
$Q(\varphi_{\delta\varepsilon}^{p-2}
\tilde{g}_{\delta\varepsilon}) \rightarrow \infty$. By the way,
the proposition \ref{bound} says that
$Y(M_{\delta\varepsilon},[\tilde{g}_{\delta\varepsilon}]_G)$ is
bounded above for any $\varepsilon\in (0,1]$, because a
$G$-invariant open set $(M_0-N(r_0),g)$ is isometrically embedded
into $(M_{\delta\varepsilon},\tilde{g}_{\delta\varepsilon})$ under
the identity map. This leads to a contradiction, completing the
proof.

\end{proof}

\begin{rmk}
A slight modification of this proof also works in more general
cases as when the two normal bundles of $W$ are isomorphic with
the equivariant $G$-action.
\end{rmk}

\section{Examples}

Consider the unit $n$-sphere $S^n(1)\subset \Bbb R^{n+1}$ for
$n\geq 3$ and an isometric $G$-action where $G=SO(n-q+1)$ with
$q\geq 3$ acts on the first $n-q+1$ coordinates of $\Bbb R^{n+1}$
fixing the last $q$ coordinates of $\Bbb R^{n+1}$. Then the
complement of fixed point set $S^n(1)\cap (\{0\}\times \Bbb R^q)$
is foliated by $G$-invariant $(n-q)$-spheres on each of which the
$G$-action is transitive. Since the round metric is a
$G$-invariant Yamabe metric and the $G$-action has fixed points,
$$Y_{G}(S^n)=Y(S^n)=\Lambda_n.$$ Take two copies of $S^n$ and
perform a surgery along such $S^{n-q}$ to get a Riemannian
$G$-manifold $S^{n-q+1}\times S^{q-1}$. By our surgery theorem,
$$Y_G(S^{n-q+1}\times S^{q-1})\geq Y_{G}(S^n\cup S^n)=
Y_{G}(S^n)=\Lambda_n.$$ Since $S^{n-q+1}\times S^{q-1}$ has fixed
points, $Y_G(S^{n-q+1}\times S^{q-1})\leq \Lambda_n$ and hence
$$Y_G(S^{n-q+1}\times S^{q-1})=\Lambda_n.$$ Taking connected sums of
$S^{n-q+1}\times S^{q-1}$ along fixed points, we also have
$$Y_G(l(S^{n-q+1}\times S^{q-1})\sharp
\ m\overline{S^{n-q+1}\times S^{q-1}})=\Lambda_n$$ for any
integers $l,m\geq 0$.


\bigskip



\begin{thebibliography}{99}

\bibitem{bes} A. Besse, Einstein Manifolds, Springer-Verlag 1987.

\bibitem{aku} K. Akutagawa and B. Botvinnik, {\em Relative Yamabe invariant},
Comm. Anal. Geom. {\bf 10} (2002) No. 5, 935--969.

\bibitem{ander} M. Anderson, {\em Remarks on Perelman's papers}, preprint.

\bibitem{aubin} T. Aubin, {\em \'{E}quations diff\'{e}rentielles
non lin\'{e}aires et probl\`{e}me de Yamabe concernant la courbure
scalaire}, J. Math. Pures Appl. {\bf 55} (1976), 269--296.

\bibitem{ber} L. B\'{e}rard Bergery, {\em Scalar curvature and isometry
groups}, in "Spectra of Riemannian manifolds", ed. M. Berger, S.
Murakami, T. Ochiai, Tokyo, (1983), 9--28.

\bibitem{Bray} H. Bray and A. Neves, {\em Prime $3$-manifolds
with Yamabe invariant greater than $\Bbb R P^3$}, Ann. of Math.
{\bf 159} (2004) No.1, 407--424.

\bibitem{GL} M. Gromov and H.B. Lawson,
 {\em  The classification of simply connected manifolds
of positive scalar curvature}, Ann. of Math. {\bf 111} (1980)
423--434.

\bibitem{HV} E. Hebey and M. Vaugon, {\em Le probl\`{e}me de
Yamabe \'{e}quivariant}, Bull. Sci. Math. {\bf 117} (1993),
241--286.

\bibitem{Hor} L. H\"{o}rmander, Linear Partial Differential
Operators, Grundlehren der math. Wiss. 116, Springer-Verlag, 1963.

\bibitem{IL} M. Ishida and C. LeBrun, {\em Curvature, connected sums,
and Seiberg-Witten theory}, Comm. Anal. Geom. {\bf 11} (2003)
No.5, 809--836.

\bibitem{joyce} D. Joyce, {\em Constant scalar curvature metrics on connected sums},
Int. J. Math. Math. Sci. (2003) No. 7, 405--450.

\bibitem{koba} O. Kobayashi, {\em Scalar curvature of a metric with unit
volume}, Math. Ann. {\bf 279} (1987), 253--265.

\bibitem{lb1} C. LeBrun,
{\em Four manifolds without Einstein metrics}, Math. Res. Lett.
{\bf 3} (1996), 133--147.

\bibitem{lb2}\leavevmode\vrule height 2pt depth -1.6pt width 23pt,
{\em Yamabe constants and the perturbed Seiberg-Witten equations},
Comm. Anal. Geom. {\bf 5} (1997), 535--553.

\bibitem{lb3}\leavevmode\vrule height 2pt depth -1.6pt width
 23pt, {\em Kodaira dimension and the Yamabe problem},
Comm. Anal. Geom.  {\bf 7} (1999), 133--156.

\bibitem{lee} J. Lee and T. Parker, {\em The Yamabe Problem},
Bull. Amer. Soc. {\bf 17} (1987), 37--81.

\bibitem{perel} G. Perelman, {\em Ricci flow with surgery on three-manifolds}, math.DG/0303109.

\bibitem{PY} J. Petean and G. Yun, {\em Surgery and the Yamabe invariant},
GAFA {\bf 9} (1999), 1189--1199.

\bibitem{sy} R. Schoen and S.T. Yau : {\em On the Structure of
Manifolds with Positive Scalar Curvature}, Manuscr. Math. {\bf 28}
(1979), 159-183.

\bibitem{Sung1} C. Sung, {\em Surgery, curvature, and minimal volume},
Ann. Global Anal. Geom. {\bf 26} (2004), 209--229.


\end{thebibliography}
\end{document}